\documentclass[12pt]{amsart}

\newtheorem{Theorem}{Theorem}[section]

\newtheorem{Corollary}[Theorem]{Corollary}
\newtheorem{Remark}[Theorem]{Remark}
\newtheorem{Example}[Theorem]{Example}
\newtheorem{Definition}[Theorem]{Definition}

\newtheorem*{ack}{\bf Acknowledgement}

\newcommand{\field}[1]{\ensuremath{\mathbb{#1}}}

\newcommand{\C}{\field{C}}

\newcommand{\N}{\field{N}}
\newcommand{\Pp}{\field{P}}
\newcommand{\Q}{\field{Q}}

\newcommand{\Z}{\field{Z}}
\newcommand{\G}{\field{G}}
\renewcommand{\P}{\field{P}}

\newcommand{\calm}{\mathcal{M}}

\newcommand{\calo}{\mathcal{O}}

\newcommand{\ps}[1]{\ensuremath{\Pp^{#1}}}

\newcommand{\tc}{{\rm TC}(X)}
\newcommand{\tcx}{\ensuremath{\bigoplus_{\overline{D} 
\in {{\rm Cl}(X)}} H^0 ( X, \calo_X (D)) }}
\newcommand{\tcxb}{\ensuremath{\bigoplus_{n_1,\ldots,n_r\in\Z}\, H^0
    \Bigl( X, \calo_X \Bigl(\sum_i n_i D_i\Bigr)\Bigr) t_1^{n_1}\cdots t_r^{n_r}}}

\DeclareMathOperator{\proj}{Proj}
\DeclareMathOperator{\cl}{Cl}

\DeclareMathOperator{\Pic}{Pic}

\DeclareMathOperator{\Hom}{Hom}

\DeclareMathOperator{\Nef}{Nef}
\DeclareMathOperator{\Ker}{Ker}

\title{The total coordinate ring of a normal projective variety}
 
\author{E. Javier Elizondo}\thanks{Partially supported by
  JSPS-CONACYT, and CONCYT 40531-F}
\address{Instituto de Matem\'aticas\\Ciudad Universitaria,
  UNAM\\M\'exico, DF 04510\\M\'exico}
\email{javier@math.unam.mx}
\author{Kazuhiko Kurano}
\address{Department of Mathematics \\Meiji University\\
Higashi-Mita~1-1-1\\ Tama-ku\\ Kawasaki 214-8571\\ Japan}
\email{kurano@math.meiji.ac.jp}
\author{Kei-ichi Watanabe}
\address{Department of Mathematics\\College of Humanities and
Sciences\\Nihon University\\Setagaya-ku\\Tokyo 156-0045\\Japan}
\email{watanabe@math.chs.nihon-u.ac.jp}
\begin{document}

\begin{flushleft}
{\it To appear in Journal of Algebra}
\end{flushleft}\vspace{.2cm}
\maketitle

\section{Introduction}
The total coordinate ring \tc\ of a variety is a generalization of
the ring introduced and studied by Cox \cite{cox-hom} in
connection with a toric variety. 
Consider a normal projective variety $X$ with divisor class group 
${\cl(X)}$,
and let us assume that it is a
finitely generated free abelian group. We define the total coordinate ring
of $X$ to be
\begin{equation*}
\tc \, = \, \ensuremath{\bigoplus_{D} H^0 ( X, \calo_X (D)) } ,
\end{equation*}
where the sum as above is taken over all Weil divisors of $X$
contained in a fixed complete system of representatives of ${\cl(X)}$.
We refer to Definition~\ref{tcr} for a precise definition of $\tc$.

Such rings grew out of an old problem of classical algebraic 
geometry, which we describe as follows.
Let $p_1$, \ldots, $p_m$ be distinct $m$ points in the projective 
space $\ps{r}$ over $\C$, and
let $\pi : X \rightarrow \ps{r}$ be the blow up of $\ps{r}$
at $\{ p_1, \ldots, p_m \}$.
Put $E_i = \pi^{-1}(p_i)$ for $i = 1, \ldots, m$.
Let $H$ be a hyperplane in $\ps{r}$ and put $A = \pi^{-1}(H)$.
Then $\cl (X)$ is a free abelian group with basis 
$\overline{E_1}$, \ldots, $\overline{E_m}$, $\overline{A}$.
Then, we may regard $H^0(X, \calo_X(\sum_{i=1}^mn_iE_i + n_{m+1}A))$  
as the linear system on $\ps{r}$ of degree $n_{m+1}$
passing through $p_i$ with multiplicity at least $-n_i$ for each $i$.
The total coordinate ring 
\[
\tc = \bigoplus_{n_1, \ldots, n_{m+1} \in \Z}
H^0\Bigl(X, \calo_X\Bigl(\sum_{i=1}^mn_iE_i + n_{m+1}A\Bigr)\Bigr)
\]
facilitates the computation of 
$H^0(X, \calo_X(\sum_{i=1}^mn_iE_i + n_{m+1}A))$, where
$\Z$ denotes the ring of integers.
In fact, if $\tc$ is a Noetherian ring, then the dimension of 
$H^0(X, \calo_X(\sum_{i=1}^mn_iE_i + n_{m+1}A))$ is a rational function
in $n_1$, \ldots, $n_{m+1}$.

In the event that $\tc$ is  Noetherian for a normal projective variety 
with function field $K$ whose divisor class group is finite free,
the ring
\[
R(D) = \bigoplus_{n \in \Z}
H^0(X, \calo_X(nD))t^n \subset K[t, t^{-1}]
\] 
is a Noetherian ring for any Weil divisor $D$ on $X$
(e.g., Elizondo and Srinivas~\cite{eli-sri}).
Furthermore,
\[
\bigoplus_{n \in \Z}
H^0(X, \calo_X(E + nD))
\] 
is a finitely generated $R(D)$-module
for any Weil divisors $D$ and $E$.

The total coordinate ring sometimes appears as an invariant subring.
Mukai exhibited numerous examples where the total coordinate rings are not
Noetherian, this being related to the  Hilbert's fourteenth problem.
See \cite{muk-14} and \cite{muk-geo}.

In the present paper, we shall prove that the total coordinate ring
 is a unique factorization domain
for any connected normal Noetherian scheme whose divisor class group is 
finitely generated free abelian group in Corollary~\ref{coro}.

The main theorem of the paper is the following:

\begin{Theorem}\label{thm1}
Let $X$ be a connected normal Noetherian scheme with function field $K$.
Let $D_1, \ldots, D_r$ be Weil divisors on $X$, and let
$t_1, \ldots, t_r$ be variables over $K$. 
We set 
\[
R = \tcxb
\subset  S = K[t_1^{\pm 1}, \ldots, t_r^{\pm 1}] .
\]
\begin{enumerate}
\item
Then, $R$ is a Krull domain.
\item
Assume that $Q(R) = Q(S)$, where $Q( \ )$ stands for the field of fractions.
Then there is a natural surjection
\[
\varphi : \cl (X) \left/ \langle \overline{D}_1, \ldots, \overline{D}_r
\rangle \right. \longrightarrow \cl (R) ,
\]
where $\overline{D}$ is the image in $\cl (X)$ of a Weil divisor $D$ on $X$.
\item
If there exist integers $n_1$, \ldots, $n_r$ such that
$\sum_in_iD_i$ is an ample divisor, then
the map $\varphi$ as above is an isomorphism.
\end{enumerate}
\end{Theorem}

The following is an immediate consequence (ref.\ Remark~\ref{rem})
of the above theorem.

\begin{Corollary}\label{coro}
With the notation as in Theorem~\ref{thm1},
$R$ is a unique factorization domain if the set of 
the images of $D_1, \ldots, D_r$ generate $\cl (X)$.

In particular, the total coordinate ring is
a unique factorization domain for a connected normal Noetherian scheme 
whose divisor class group is a finitely generated free abelian 
group.
\end{Corollary}

We should mention that independently
Berchtold and Hausen~\cite{BH} proved that
the total coordinate ring is a unique factorization domain
for a locally factorial variety over an algebraically closed field
whose Picard group is finitely generated
free abelian group.
Our method is purely algebraic, simple, and totally different from theirs.
Here, we give a precise statement of the result of 
Berchtold and Hausen~\cite{BH}.
Let $X$ be a normal variety over an algebraically closed field 
whose Picard group is finitely generated free abelian group.
Let
\[
{\rm A}(X) = \bigoplus_{\overline{D} \in {\rm Pic}(X)}
H^0(X, \calo_X(D))
\]
be the ring defined in the same way as in Definition~\ref{tcr} using 
all Cartier divisors.
Observe that if the divisor class group of $X$ is finitely generated
free abelian group, this is the subring
of ${\rm TC}(X)$ [see (\ref{tc}) in Definition~\ref{tcr}] 
consisting of those graded
pieces corresponding to Cartier divisors.
Then, Proposition~8.4 in \cite{BH}
implies that the following four conditions are equivalent;
(1) ${\rm A}(X)$ is a unique factorization domain,
(2) $X$ is locally factorial,
(3) ${\rm Cl}(X) = {\rm Pic}(X)$,
(4) ${\rm A}(X) = {\rm TC}(X)$.
Let $X_{\rm reg}$ be the open subscheme consisting of non-singular points
of $X$.
Furthermore, assume that the divisor class group of $X$ is finitely generated
free abelian group.
Since the codimension of $X \setminus X_{\rm reg}$ in $X$ is at least $2$,
${\rm TC}(X) = {\rm TC}(X_{\rm reg})$ is satisfied.
Then, Proposition~8.4 in \cite{BH} implies that
${\rm TC}(X_{\rm reg})$ is a unique
factorization domain since $X_{\rm reg}$ is locally factorial.
Therefore we know that 
${\rm TC}(X)$ is a unique factorization domain by the result of
Berchtold and Hausen.
Here, we remark that ${\rm A}(X_{\rm reg})$ coincides with
${\rm TC}(X_{\rm reg})$ since $X_{\rm reg}$ is locally factorial.
On the other hand, ${\rm A}(X) = {\rm TC}(X)$ if and only if
$X$ is locally factorial.

We remark that we do not assume that a scheme is of finite type over
an algebraically closed field in Corollary~\ref{coro}.

We shall prove Theorem~\ref{thm1} in the next section.
In the final section, we give some examples of total coordinate rings.
We prove that the total coordinate ring of the blow up 
of a projective space at a finite number of points on a line is finitely
generated, see Example~\ref{noet}.
Next we give an example of a normal projective variety $X$
that has a non-finitely generated total coordinate ring
in Example~\ref{example}.

\begin{ack}
\begin{rm}
The first author would like to thank the hospitality 
of the Tokyo Metropolitan University, during his visit in October 2001,  
where the proof of the main theorems
were discussed. The second and third authors visited Mexico in 2000.
We would like
to thank the exchange program between Mexico and Japan,
JSPS-Conacyt which supported our visit. We also would like to
thank James Lewis for many corrections, to Paul Roberts and
Vasudevan Srinivas for reading the first draft of this article, and
to Florian Berchtold, J\"urgen Hausen and the referee for
many valuable comments.
\end{rm}
\end{ack}

\section{The total coordinate ring}
\label{sct1}

In this section we define the total coordinate ring 
and prove Theorem~\ref{thm1}.

\begin{Definition}
\label{tcr}
\begin{rm}
Let $X$ be a connected normal Noetherian scheme with function field $K$. We
assume that the divisor class
group ${\cl (X)}$ of $X$ is isomorphic to $\Z^{\oplus s}$. 
Let $D_1$, \ldots, $D_s$ be Weil divisors on $X$ such that the set of
their images generate ${\cl (X)}$.
Let $t_1$, \ldots, $t_s$ be variables over $K$.
We set 
\[
R(X;D_1, \ldots, D_s) = \bigoplus_{n_1, \ldots, n_s \in \Z}
H^0\Bigl(X, \calo_X \Bigl(\sum_in_iD_i\Bigr)\Bigr) t_1^{n_1} \cdots t_s^{n_s}
\subset K[t_1^{\pm 1}, \ldots, t_s^{\pm 1}] ,
\]
where we regard
\[
H^0\Bigl(X, \calo_X \Bigr(\sum_in_iD_i\Bigr)\Bigr) = 
\{ a \in K^\times \mid {\rm div}_X(a) + \sum_in_iD_i \geq 0 \} \cup \{ 0 \}
\]
as an additive subgroup of $K$.
It is easily seen that $R(X;D_1, \ldots, D_s)$ is a subring of
$K[t_1^{\pm 1}, \ldots, t_s^{\pm 1}]$.
The ring $R(X;D_1, \ldots, D_s)$ is 
uniquely determined
by $X$ up to isomorphism, that is, it is independent of the choice of
$D_1$, \ldots, $D_s$ up to isomorphism.
We call it the {\it total coordinate ring}  $\tc$ of $X$.
We may regard the total coordinate ring of $X$ as a graded ring
with the grading given by ${\cl (X)}$.
We sometimes write it as
\begin{equation}
\label{tc}
\tc \, = \, \tcx ,
\end{equation}
where $\overline{D}$ denotes the class in ${\cl (X)}$ represented 
by a Weil divisor $D$.
\end{rm}
\end{Definition}

Before proving the main theorem, we give a remark.

\begin{Remark}
\label{rem}
\begin{rm}
With the notation given in Theorem~\ref{thm1},
it is easily seen that $Q(R)$ coincides with $Q(S)$
if one of the following two conditions are satisfied;
(1) the set of the images of $D_1$, \ldots, $D_s$ generate ${\rm Cl}(X)$,
(2) there exist integers $n_1$, \ldots, $n_r$ such that
$\sum_in_iD_i$ is an ample divisor.

Even if $Q(R) = Q(S)$ is satisfied, the map $\varphi$ in Theorem~\ref{thm1} 
is not an isomorphism
in general, as can
be seen in the following example.
Let $\pi : X \rightarrow \ps{2}$ be the blow up of $\ps{2}$
at a point $p$.
Put $E = \pi^{-1}(p)$.
Let $H$ be a hyperplane in $\ps{2}$ and put $A = \pi^{-1}(H)$.
Let $K$ be the function field of $X$ and let $t$ be a variable over $K$.
Set 
\[
R = \bigoplus_{n \in \Z}
H^0(X, \calo_X (nA)) t^n
\subset S = K[t^{\pm 1}] .
\]
Using the projection formula, we have 
$H^0(X, \calo_X (nA)) = H^0(\ps{2}, \calo_{\ps{2}} (nH))$.
Therefore $Q(R) = Q(S)$ is satisfied.
In this case, $\cl(X)$ is a $\Z$-free module of rank $2$ with basis 
$\overline{E}$ and $\overline{A}$.
Since $R$ is a polynomial ring over the base field, 
we have ${\rm Cl}(R) = 0$.
Therefore the map $\varphi$ in Theorem~\ref{thm1} is not an isomorphism
in this case.
\end{rm}
\end{Remark}

For the remainder of this section, we prove Theorem~\ref{thm1}.

Let $H_1$ be the set of reduced and irreducible closed subschemes
of $X$ of codimension $1$.
Put
\[
D_i = \sum_{F \in H_1} m_{i,F} F
\]
for $i = 1, \ldots, r$.
Let $\calo_{X,F}$ be the local ring of $X$ at $F$, and let $\calm_{X,F}$ be
the corresponding maximal ideal.
For $F \in H_1$, we denote by $v_F$ the normalized valuation of the 
discrete valuation ring $\calo_{X,F}$.
Then we have
\[
R = \bigoplus_{n_1, \ldots, n_r \in \Z}
\{ a \in K \mid 
\mbox{$v_F(a) + \sum_in_im_{i,F} \geq 0$ for each $F \in H_1$} \}
t_1^{n_1} \cdots t_r^{n_r} .
\]
Here we set
\[
R_F = \bigoplus_{n_1, \ldots, n_r \in \Z}
\{ a \in K \mid 
\mbox{$v_F(a) + \sum_in_im_{i,F} \geq 0$} \}
t_1^{n_1} \cdots t_r^{n_r} 
\]
for $F \in H_1$.
It is easy to see that $R_F$ is also a subring of 
$S = K[t_1^{\pm 1}, \ldots, t_r^{\pm 1}]$ such that
$R \subset R_F \subset S$ for any $F \in H_1$,
and that $R = \cap_{F \in H_1}R_F$ is satisfied.

Since the local ring $(\calo_{X,F}, \calm_{X,F})$ is 
a discrete valuation ring,
there exists an element $\alpha_F \in \calm_{X,F}$ such that
$\calm_{X,F} = \alpha_F \calo_{X,F}$.
By the equality
\[
\{ a \in K \mid 
\mbox{$v_F(a) + \sum_in_im_{i,F} \geq 0$} \}
= \alpha_F^{-\sum_in_im_{i,F}}\calo_{X,F} ,
\]
we have
\begin{eqnarray*}
R_F & = & \bigoplus_{n_1, \ldots, n_r \in \Z}
\alpha_F^{-\sum_in_im_{i,F}}\calo_{X,F} t_1^{n_1} \cdots t_r^{n_r} \\
& = & \bigoplus_{n_1, \ldots, n_r \in \Z}
\calo_{X,F} (\alpha_F^{-m_{1,F}}t_1)^{n_1} \cdots 
(\alpha_F^{-m_{r,F}}t_r)^{n_r} \\
& = & \calo_{X,F}\left[(\alpha_F^{-m_{1,F}}t_1)^{\pm 1}, \ldots,
(\alpha_F^{-m_{r,F}}t_r)^{\pm 1}
\right] .
\end{eqnarray*}
Therefore, $R_F$ is a Noetherian normal domain.
We remark that $\alpha_F R_F$ is the unique homogeneous prime ideal
of height $1$ of $R_F$.
Let $\{ Q_\lambda \mid \lambda \in \Lambda \}$ be the set of non-homogeneous
prime ideals of height $1$ of $R_F$.
Since $R_F$ is a Krull domain,
\[
R_F = (R_F)_{\alpha_F R_F} \bigcap \left( \bigcap_{\lambda \in \Lambda}
(R_F)_{Q_\lambda} \right)
\]
is satisfied (see Theorem~12.3 in Matsumura \cite{mat-com}).
Since $S = R_F[\alpha_F^{-1}]$, we have
\[
S = \bigcap_{\lambda \in \Lambda} (R_F)_{Q_\lambda} 
\]
and
\begin{equation}
\label{eq1}
R = \cap_{F \in H_1}R_F
= \left( \bigcap_{F \in H_1}
(R_F)_{\alpha_F R_F} \right) \cap S .
\end{equation}

Here for $F \in H_1$, we set $P_F = R \cap \alpha_F R_F$.
Then we have
\begin{eqnarray}
\nonumber
P_F & = & \bigoplus_{n_1, \ldots, n_r \in \Z}
\left\{ a \in K \left|
\begin{array}{ll}
\mbox{$v_F(a) + \sum_in_im_{i,F} > 0$, and} \\
\label{P_F}
\mbox{$v_G(a) + \sum_in_im_{i,G} \geq 0$ for each $G \in H_1$}
\end{array}
\right.
\right\}
t_1^{n_1} \cdots t_r^{n_r} \\
& = & \bigoplus_{n_1, \ldots, n_r \in \Z}
H^0(X, \calo_X(\sum_in_iD_i - F)) t_1^{n_1} \cdots t_r^{n_r} .
\end{eqnarray}

We shall prove that $R$ is a Krull domain.
By equation (\ref{eq1}), it is enough to show that
for any $a \in R \setminus \{ 0 \}$, 
$a$ is a unit in $(R_F)_{\alpha_F R_F}$
except for finitely many $F$'s.
(Here, we remark that $Q(R)$ does not have to coincides with $Q(S)$.)
We note that
\[
\alpha_F (R_F)_{\alpha_F R_F} \cap R 
= \alpha_F R_F \cap R = P_F .
\]
By equation~(\ref{P_F}) as above, it is easy to see that there exist
only finitely many $F$'s such that $a$ is contained in $P_F$.
We have thus proven that $R$ is a Krull domain.

%

By Theorem~12.3 in Matsumura \cite{mat-com}, 
the set $\{ P_F \mid F \in H_1 \}$ includes the set of homogeneous prime ideals
of height $1$ of $R$.
Set
\[
H'_1 = \{ F \in H_1 \mid {\rm ht}_R P_F = 1 \} .
\]
Thus $\{ P_F \mid F \in H'_1 \}$ is the set of homogeneous prime ideals 
of height $1$ of $R$.
(Let $E$ be the divisor given in Remark~\ref{rem}.
Then, it is easy to see that ${\rm ht}_RP_E = 2$.
Therefore $E \in H_1 \setminus H'_1$ in this case.)
By Theorem~12.3 in \cite{mat-com},
we have
\[
R 
= \left( \bigcap_{F \in H'_1}
R_{P_F} \right) \cap S .
\]

We denote by ${\rm Div}(X)$ the group of Weil divisors on $X$.
It is the free abelian group generated by $H_1$.
Let ${\rm P}(X)\subset {\rm Div}(X)$ be
the subgroup of principal divisors
\[
\{ {\rm div}_X(a) \mid a \in K^{\times} \} ,
\]
where
\[
{\rm div}_X(a) = \sum_{F \in H_1} v_F(a) F .
\]
By definition, $\cl(X)$ is the quotient group ${\rm Div}(X)/{\rm P}(X)$.
Let $v'_F$ be the normalized valuation of the discrete valuation ring
$R_{P_F}$ for $F \in H'_1$.

For the remainder of the proof, 
we assume that $Q(R)$ coincides with $Q(S)$.
We remark that $R_{P_F} = (R_F)_{\alpha_F R_F}$
is satisfied for $F \in H'_1$ since $Q(R) = Q(S)$.
We denote by ${\rm HDiv}(R)$ the set of homogeneous Weil divisors on $R$. 
That is, it is the free abelian group generated by $\{ P_F \mid F \in H'_1 \}$.
Let ${\rm HP}(R)\subset {\rm HDiv}(R)$ denote the subgroup 
generated by 
\[
\{ {\rm div}_R(b) \mid \mbox{$b$ is a non-zero homogeneous 
element of $R$} \} ,
\]
where
\[
{\rm div}_R(b) = \sum_{F \in H'_1} v'_F(b) P_F .
\]
Then, it is known that the natural map 
\[
{\rm HDiv}(R)/{\rm HP}(R) \rightarrow \cl(R)
\]
is an isomorphism (e.g., see Samuel~\cite{Sam64}).

Let $\xi : {\rm Div}(X) \rightarrow {\rm HDiv}(R)$ be the surjective
homomorphism 
given by $\xi(F) = P_F$ for $F \in H'_1$ and $\xi(F) = 0$ for
$F \in H_1 \setminus H'_1$.
Consider the following diagram:
\[
\begin{array}{ccccccccc}
0 & \longrightarrow & {\rm P}(X) & \longrightarrow & {\rm Div}(X) & 
\longrightarrow & \cl(X) & \longrightarrow 0 \\
& & & & \scriptstyle{\xi} \downarrow \phantom{\scriptstyle{\xi}}
& & & & \\
0 & \longrightarrow & {\rm HP}(R) & \longrightarrow & {\rm HDiv}(R) & 
\longrightarrow & \cl(R) & \longrightarrow 0 \\
\end{array}
\]
For $F \in H'_1$ and $a \in K^{\times}$,
$v_F(a) = v'_F(a)$ is satisfied since $R_{P_F} = (R_F)_{\alpha_F R_F}$.
Therefore for $a \in K^{\times}$,
\[
\xi({\rm div}_X(a)) = {\rm div}_R(a)
\]
is satisfied.
Hence we have $\xi({\rm P}(X)) \subset {\rm HP}(R)$.

By definition, the additive group
${\rm HP}(R)$ is generated by
\[
\{ {\rm div}_R(a) \mid a \in K^{\times} \} 
\mbox{ \ \ and \ \ }
\{ {\rm div}_R(t_i) \mid i = 1, \ldots, r \} ,
\]
since $Q(R) = Q(S)$.
For $F \in H'_1$,
since $R_{P_F} = (R_F)_{\alpha_F R_F}$ and
\[
R_F = 
\calo_{X,F}\left[(\alpha_F^{-m_{1,F}}t_1)^{\pm 1}, \ldots,
(\alpha_F^{-m_{r,F}}t_r)^{\pm 1}
\right] ,
\]
$v'_F(\alpha_F^{-m_{i,F}}t_i) = 0$ is satisfied.
Therefore we have $v'_F(t_i) = m_{i,F}$.
Hence for each $i$, 
\[
\xi(D_i) = \xi(\sum_{F \in H_1} m_{i,F} F)
= \sum_{F \in H'_1} v'_F(t_i) P_F = {\rm div}_R(t_i) 
\]
is satisfied.
Since ${\rm HP}(R)$ is generated by $\xi({\rm P}(X))$ and
$\{ \xi(D_i) \mid i = 1, \ldots, r \}$,
we have an isomorphism
\begin{eqnarray*}
\cl(R) & \simeq & 
{\rm Div}(X)\left/
{\rm P}(X) + \langle D_1, \ldots, D_r\rangle +
\langle F \mid F \in H_1 \setminus H'_1 \rangle
\right. \\
& = &
\cl(X) \left/
\langle \overline{D_1}, \ldots, \overline{D_r}\rangle +
\langle \overline{F} \mid F \in H_1 \setminus H'_1 \rangle
\right.   .
\end{eqnarray*}
By this isomorphism, we arrive at the natural surjection
\[
\varphi : \cl (X) \left/ \langle \overline{D_1}, \ldots, \overline{D_r}
\rangle \right. \longrightarrow \cl (R) .
\]

For the rest of this section, we shall prove that $\varphi$ is an isomorphism
if the subgroup $\langle D_1, \ldots, D_r \rangle$ of ${\rm Div}(X)$
contains an ample divisor.
We remark that $Q(R) = Q(S)$ is satisfied in this case.
It is enough to show $H_1 = H'_1$.
In order to show this, we need only  show that 
$R_{P_F} = (R_F)_{\alpha_F R_F}$
for each $F \in H_1$.
It is sufficient to show $R_F \subset R_{P_F}$ for each $F \in H_1$.

Let $f$ be a homogeneous element of $R_F$.
Since $R_F$ is a graded ring such that $R \subset R_F \subset Q(R)$,
it is easy to see that there exist $\alpha, \beta \in K^{\times}$ and integers
$p_1$, \ldots, $p_r$, $q_1$, \ldots, $q_r$ such that
\begin{equation*}
\label{eq3}
\alpha t_1^{p_1}\cdots t_r^{p_r}, \beta t_1^{q_1}\cdots t_r^{q_r} \in R
\ \ \mbox{and} \ \ f = 
\frac{\beta t_1^{q_1}\cdots t_r^{q_r}}{\alpha t_1^{p_1}\cdots t_r^{p_r}} 
\in R_F .
\end{equation*}
Thus we have
\begin{eqnarray*}
{\rm div}_X(\alpha) + \sum_ip_iD_i & \geq & 0 \\
{\rm div}_X(\beta) + \sum_iq_iD_i & \geq & 0 \\
v_F(\beta/\alpha) + \sum_i(q_i - p_i)m_{i,F} & \geq & 0 .
\end{eqnarray*}
We want to show $f \in R_{P_F}$.

If $v_F(\alpha) + \sum_ip_im_{i,F} = 0$, then 
$\alpha t_1^{p_1}\cdots t_r^{p_r} \in R
\setminus P_F$ is satisfied.
Therefore in this case, we have $f \in R_{P_F}$.

Assume that $v_F(\alpha) + \sum_ip_im_{i,F} > 0$.
Put $p = v_F(\alpha) + \sum_ip_im_{i,F}$.
We remark that $v_F(\beta) + \sum_iq_im_{i,F} \geq p$.
Since $\langle D_1, \ldots, D_r \rangle$ contains an ample divisor,
there exist $\gamma \in K^{\times}$ and integers $s_1$, \ldots, $s_r$ such that
\begin{eqnarray*}
{\rm div}_X(\gamma) + \sum_is_iD_i + pF & \geq & 0 \\
v_F(\gamma) + \sum_is_im_{i,F} + p & = & 0 .
\end{eqnarray*}
Then we have
\begin{eqnarray*}
{\rm div}_X(\beta\gamma) + \sum_i(q_i + s_i)D_i & \geq & 0 \\
{\rm div}_X(\alpha\gamma) + \sum_i(p_i + s_i)D_i & \geq & 0 \\
v_F(\alpha\gamma) + \sum_i(p_i + s_i)m_{i,F} & = & 0 .
\end{eqnarray*}
Therefore we have 
$\beta\gamma t_1^{q_1+s_1}\cdots t_r^{q_r+s_r} \in R$
and
$\alpha\gamma
t_1^{p_1+s_1}\cdots t_r^{p_r+s_r} \in R \setminus P_F$.
Since
\[
f  = 
\frac{\beta t_1^{q_1}\cdots t_r^{q_r}}{\alpha t_1^{p_1}\cdots t_r^{p_r}} =
\frac{\beta\gamma t_1^{q_1+s_1}\cdots t_r^{q_r+s_r}}{\alpha\gamma
t_1^{p_1+s_1}\cdots t_r^{p_r+s_r}} ,
\]
we have $f \in R_{P_F}$.

This completes the proof of Theorem~\ref{thm1}.

\begin{Remark}
\begin{rm}
Yi Hu and Se\'an Keel~\cite[Theorem 2.9]{kee-git} proved the following.

 Let $X$ be a
\Q-factorial projective variety such that $\Pic(X)_{\Q} \, =\, {\rm
  N}^1(X)$. Then $X$ is a Mori dream space if and only if \tc\ is
finitely generated. If $X$ is a Mori dream space then
$X$ is a GIT quotient of $V\, = \, spec(\tc) $ by the torus $G\, = \,
\Hom (\N^r, \G_m)$

Here a Mori dream space is a variety with nice geometric properties.
For example, the nef cone
$\Nef(X)$ is the affine hull of finitely many semi-ample line bundles,
and there exist small \Q-factorial modifications of $X$.
\end{rm}
\end{Remark}

\section{Some examples}

We give some examples of total coordinate rings in the section.

\begin{Example}
\begin{rm}
It is well known that the divisor class group is a finitely generated 
free abelian group for a smooth complete toric variety 
(e.g., see 63p in \cite{ful-tova}).
Furthermore in this case,
Cox~\cite{cox-hom} proved that
\tc\ is a homogeneous polynomial ring.  
He called \tc\ the {\it homogeneous
coordinate ring} of $X$. 
\end{rm}
\end{Example}

\begin{Remark}
\begin{rm}
Total coordinate rings have a deep relation with 
invariant theory as Mukai \cite{muk-14} has shown. Here we present some
 of his results.

Let $\ps{r}$ be the projective space of dimension $r$ over the field
of complex numbers ${\C}$.
Let $X$ be the blow up of $\ps{r}$ at
$m$ distinct points.
Then we have ${\cl}(X) \simeq {\Z}^{m+1}$. 
Assume that the $m$ points are not on any hyperplane in $\ps{r}$.
Then with a suitable linear action of $G = {\mathbf G}_a^{m-r-1}$ on 
a polynomial ring $S$ over ${\C}$ with $2m$ variables,
the invariant subring $S^G$ is isomorphic to ${\tc}$. 

On the other hand,
Nagata \cite{nag-rat} proved that 
the effective cone ${\rm Eff}(X)$ in $\cl(X)$ is not finitely generated as 
a semi-group if $X$ is a blow up of $\ps{2}$ at
$9$ general points.
Therefore, $\tc$ is not Noetherian in this case.
It is a counterexample of Hilbert's 14-th problem. 

Mukai also  generalized in \cite{muk-geo} a result of 
Dolgachev and realized the
root system $T_{p,q,r}$ in the cohomology group of a certain rational
variety of Picard number $p+q+r-1$. As an application he proved that
the invariant subring of a tensor product with an actions of Nagata type
is infinitely generated if the Weyl group of the corresponding root
system $T_{p,q,r}$ is infinite. 
\end{rm}
\end{Remark}

\begin{Example}\label{noet}
\begin{rm}
Let $p_1, \ldots, p_m$ be $m$ points on a projective line contained in
the projective $r$-space
\ps{r} over an algebraically closed field $k$. Then the
 total coordinate ring \tc\ of the blow up $X=
Bl_{p_1,\ldots,p_m}(\ps{r})$ of \ps{r}\ at $\{ p_1, \ldots,
p_m\}$ is a Noetherian ring.   

We give a proof:

Let $p_1$, \ldots, $p_m$ be $m$ distinct points in $\ps{r}$.
Let $\pi : X \rightarrow \ps{r}$ be the blow up of $\ps{r}$
at $\{ p_1, \ldots, p_m \}$.
Put $E_i = \pi^{-1}(p_i)$ for $i = 1, \ldots, m$.
Let $H$ be a hyperplane in $\ps{r}$ and put $A = \pi^{-1}(H)$.
Then $\cl (X)$ is a free abelian group with basis 
$\overline{E_1}$, \ldots, $\overline{E_m}$, $\overline{A}$.

Let $B = {k}[Z_0, \ldots, Z_r]$ be the homogeneous coordinate ring
of $\ps{r}$.
We denote by $I_i$ the homogeneous prime ideal of $B$ corresponding to 
the point $p_i$.
For $i = 1, \ldots, m$ and $s \in {\Z}$, we define
\[
F_{is} = \left\{
\begin{array}{ll}
I_i^s & \mbox{if $s \geq 0$} \\
B & \mbox{if $s < 0.$}
\end{array}
\right. 
\]

Then we have
\[
[F_{1b_1} \cap \cdots  \cap F_{mb_m}]_a = 
H^0(X, {\calo}_X(aA - b_1E_1 - \cdots - b_mE_m)) ,
\]
where $[ \ ]_a$ denotes the homogeneous component of degree $a$.
Therefore we have 
\begin{eqnarray*}
\tc & = & \bigoplus_{a, b_1, \ldots, b_m \in {\Z}}
H^0(X, {\calo}_X(aA - b_1E_1 - \cdots - b_mE_m)) \\
& = & \bigoplus_{b_1, \ldots, b_m \in {\Z}}
(F_{1b_1} \cap \cdots  \cap F_{mb_m})T_1^{b_1} \cdots T_m^{b_m} \\
& \subset & B[T_1^{\pm 1}, \ldots, T_m^{\pm m}]
\end{eqnarray*}

For the rest of the proof, we assume that $\{ p_1, \ldots, p_m \}$
lie on a line in $\ps{r}$. 
We may assume 
\begin{eqnarray*}
I_1 & = & (f_1, Z_2, \ldots, Z_r) \\
I_2 & = & (f_2, Z_2, \ldots, Z_r) \\
& \vdots & \\
I_m & = & (f_m, Z_2, \ldots, Z_r) , \\
\end{eqnarray*}
where $f_1, \cdots, f_m \in {k}[Z_0, Z_1]$ are linear forms
such that any two elements in $\{ f_1, \cdots, f_m \}$ are linearly independent
over ${k}$.

Here, for $i = 1, \ldots, m$ and $s \in {\Z}$, we put
\[
f_{is} = \left\{
\begin{array}{ll}
f_i^s & \mbox{if $s \geq 0$} \\
1 & \mbox{if $s < 0.$}
\end{array}
\right. 
\]

Let $g$ be a polynomial in $B$.
Put 
\[
g = \sum_{c_2, \ldots, c_r \geq 0}g_{c_2 \cdots c_r}Z_2^{c_2}\cdots Z_r^{c_r} ,
\]
where $g_{c_2 \cdots c_r} \in {k}[Z_0, Z_1]$.
Then $g$ is contained in $F_{ib_i}$ if and only if  
$g_{c_2 \cdots c_r} \in f_{i, b_i - c_2 - \cdots - c_r}{k}[Z_0, Z_1]$
for any integers $c_2, \ldots, c_r \geq 0$.
Here, $b_i$ or $b_i - c_2 - \cdots - c_r$ are possibly negative.
Then $g$ is contained in $F_{1b_1} \cap \cdots  \cap F_{mb_m}$
if and only if
 $g_{c_2 \cdots c_r} \in f_{i, b_i - c_2 - \cdots - c_r}{k}[Z_0, Z_1]$
for any integers $c_2, \ldots, c_r \geq 0$ and $i = 1, \ldots, m$.
Furthermore, this is equivalent to 
\[
g_{c_2 \cdots c_r} \in f_{1, b_1 - c_2 - \cdots - c_r} \cdots
f_{m, b_m - c_2 - \cdots - c_r}{k}[Z_0, Z_1]
\]
for any integers $c_2, \ldots, c_r \geq 0$.
In particular, $g \in F_{1b_1} \cap \cdots  \cap F_{mb_m}$
if and only if $g_{c_2 \cdots c_r}Z_2^{c_2}\cdots Z_r^{c_r} \in 
F_{1b_1} \cap \cdots \cap F_{mb_m}$ for any $c_2, \ldots, c_r \geq 0$.

Here we claim that
\begin{equation}
\label{tousiki}
\tc = B[T_1^{-1}, \ldots, T_m^{-1}, Z_2T_1 \cdots T_m,
\ldots, Z_rT_1 \cdots T_m, f_1T_1, \ldots, f_mT_m ] .
\end{equation}
It is easy to see that the right-hand side is included in the left 
side.
Assume that $g \in F_{1b_1} \cap \cdots  \cap F_{mb_m}$.
Then $gT_1^{b_1} \cdots T_m^{b_m}$ is contained in the left side.
We want to show that the right side also contains it.
We may assume that $g = g_1Z_2^{c_2}\cdots Z_r^{c_r}$,
where $g_1 \in {k}[Z_0, Z_1]$ and $c_1, \ldots, c_r$ are non-negative 
integers.
Since $g_1Z_2^{c_2}\cdots Z_r^{c_r} \in F_{1b_1} \cap \cdots  \cap F_{mb_m}$,
we may assume that 
\[
g_1 = g_2 f_{1, b_1 - c_2 - \cdots - c_r} \cdots
f_{m, b_m - c_2 - \cdots - c_r} ,
\]
where $g_2 \in {k}[Z_0, Z_1]$.
Then
\begin{eqnarray*}
& & gT_1^{b_1} \cdots T_m^{b_m} \\
& = & g_2 f_{1, b_1 - c_2 - \cdots - c_r} \cdots
f_{m, b_m - c_2 - \cdots - c_r}Z_2^{c_2}\cdots Z_r^{c_r}
T_1^{b_1} \cdots T_m^{b_m} \\
& = & g_2 (Z_2T_1 \cdots T_m)^{c_2} \cdots (Z_rT_1 \cdots T_m)^{c_r} \\
& & \times 
(f_{1, b_1 - c_2 - \cdots - c_r}T_1^{b_1 - c_2 - \cdots - c_r})
\cdots
(f_{m, b_m - c_2 - \cdots - c_r}T_m^{b_m - c_2 - \cdots - c_r}) .
\end{eqnarray*}
Here, if $b_i - c_2 - \cdots - c_r < 0$, then
\[
f_{i, b_i - c_2 - \cdots - c_r}T_i^{b_i - c_2 - \cdots - c_r}
= (T_i^{-1})^{c_2 + \cdots + c_r - b_i} .
\]
If $b_i - c_2 - \cdots - c_r \geq 0$, then
\[
f_{i, b_i - c_2 - \cdots - c_r}T_i^{b_i - c_2 - \cdots - c_r}
= (f_iT_i)^{b_i - c_2 - \cdots - c_r} .
\]
Thus $gT_1^{b_1} \cdots T_m^{b_m}$ is contained
in the ring on the right-hand side in (\ref{tousiki}), and 
this completes the proof.

Observe that if $m \geq 2$, then we obtain
\begin{eqnarray*}
& & B[T_1^{-1}, \ldots, T_m^{-1}, Z_2T_1 \cdots T_m,
\ldots, Z_rT_1 \cdots T_m, f_1T_1, \ldots, f_mT_m ] \\
& = & {k}[T_1^{-1}, \ldots, T_m^{-1}, Z_2T_1 \cdots T_m,
\ldots, Z_rT_1 \cdots T_m, f_1T_1, \ldots, f_mT_m ] .
\end{eqnarray*}
\end{rm}
\end{Example}

We shall give an example of a normal projective variety 
with infinitely generated total coordinate ring.

\begin{Example}
\label{example}
\begin{rm}
Let us consider the weighted polynomial ring $B_k$ over a field $k$
in three variables $x,y,z$
of degree $a,b,c$, respectively.
Take the weighted projective plane
\[
\P^2_k(a,b,c) = \proj (B_k)
\] 
and consider the blow up  $\pi : X_k(a,b,c) \rightarrow \P^2_k(a,b,c)$ 
at the smooth point
\[
{\mathfrak p}_k(a,b,c) := 
\Ker \bigl(k[x,y,z] \stackrel{\varphi}{\longrightarrow} k[t]\bigr) ,
\]
where $\varphi$ is the homomorphism of $k$-algebras
defined by $\varphi(x) = t^a$, $\varphi(y) = t^b$ and $\varphi(z) = t^c$.
We denote $X_k(a,b,c)$ and ${\mathfrak p}_k(a,b,c)$
simply by $X_k$ and ${\mathfrak p}_k$, respectively.
Using intersection theory on $X_k$,
Cutkosky~\cite{cut} studied the finite generation of
the symbolic Rees ring $R_s({\mathfrak p}_k) = 
\oplus_{n \geq 0}{\mathfrak p}_k^{(n)}T^n \subset B_k[T]$.
We remark that the total coordinate ring 
${\rm TC}(X_k)$ is equal to $R_s({\mathfrak p}_k)[T^{-1}]$.
Furthermore, ${\rm TC} (X_k)$ is Noetherian ring
if and only if $R_s({\mathfrak p}_k)$ is.

Assume that 
$a=7n-3,b=n(5n-2),c=8n-3$ with $n\ge 4$ and $(n,3)=1$.
In this case, the total coordinate ring ${\rm TC} (X_k)$
is a Noetherian ring if and only if the characteristic of $k$ is positive 
(see Goto, Nishida and Watanabe~\cite{GNW94}).
\end{rm}
\end{Example}


\end{document}